\documentclass{article}
\usepackage{graphicx,color}
\usepackage{amssymb}
\usepackage[geometry]{ifsym}

\newtheorem{lemma}{Lemma}
\newtheorem{theorem}{Theorem}

\newtheorem{conjecture}{Conjecture}
\newcommand{\DONOTTEX}[1]{}

\begin{document}

\title{{\bf On minimal graphs \\ containing \boldmath$k$ perfect matchings}}
\author{Ga\v{s}per Fijav\v{z} \and Matthias Kriesell}

\maketitle

\begin{abstract}
  \setlength{\parindent}{0em}
  \setlength{\parskip}{1.5ex}

  We call a finite undirected graph {\em minimally $k$-matchable} if it has at least
  $k$ distinct perfect matchings but deleting any edge results in a graph which has not.
  An {\em odd subdivision} of some graph $G$ is any graph obtained by replacing
  every edge of $G$ by a path of odd length connecting its endvertices
  such that all these paths are internally disjoint. 
  We prove that for every $k \geq 1$ there exists a finite set of graphs ${\frak G}_k$ such
  that every minimally $k$-matchable graph is isomorphic to a disjoint union of an odd subdivision of some graph from ${\frak G}_k$ and any number of copies of $K_2$.
          

  
\end{abstract}

\maketitle

\section{Introduction}

All {\em graphs} considered here are supposed to be finite and undirected unless stated otherwise, and they may contain multiple edges but no loops.
For terminology not defined here, we refer to \cite{BondyMurty2007} or \cite{Diestel2010}.
A {\em matching} of $G$ is a set $M$ of edges of $G$ such that every vertex of $G$ is end vertex of at most one member of $M$,
and $M$ is called a {\em perfect} matching of $G$ if every vertex of $G$ is end vertex of exactly one member of $M$.
By ${\frak M}(G)$ we denote the set of perfect matchings of $G$.
A graph is {\em $k$-matchable} if $|{\frak M}(G)| \geq k$, and it is called 
{\em minimally $k$-matchable} if it is $k$-matchable but, for every $e \in E(G)$, $G-e$ is not.
An {\em odd subdivision} (sometimes called a {\em totally odd subdivision}) of a graph $G$ is any graph obtained from $G$ by replacing every $e$ edge with a path of odd length (possibly $1$)
connecting the end vertices of $e$ such that all these paths are pairwise internally disjoint.
In particular, $G$ is an odd subdivision of itself. Our main result is the following.

\begin{theorem}
  \label{T1}
  For every $k \geq 1$ there exists a finite set of graphs ${\frak G}_k$ such that every minimally $k$-matchable graph
  is isomorphic to the disjoint union of an odd subdivision of some graph from ${\frak G}_k$ and any number of copies of $K_2$.
\end{theorem}

It is easy to see that the minimally $1$-matchable graphs are just disjoint unions of any number (perhaps $0$) of copies of $K_2$, and that
the minimally $2$-matchable graphs are disjoint unions of a single cycle $C_\ell$ of even length $\ell \geq 2$ and any number of copies of $K_2$.
So Theorem \ref{T1} holds with ${\frak G}_1=\emptyset$, and for ${\frak G}_2=\{C_2\}$.
However, the situation gets more complex for larger $k$, not only in terms of an increasing size of the sets ${\frak G}_k$;
for example, the classes of minimally $k$-matchable graphs need not even to be disjoint for distinct $k$:
The disjoint union $G$ of two even cycles has four perfect matchings, but deleting any edge results in a graph which has only two perfect matchings;
therefore, $G$ is minimally $4$-matchable and, at the same time, minimally $3$-matchable.

There are some results on graphs with a fixed number of perfect matchings.
For example it is known that for every positive integer $k$ there exists a constant $c_k$ such that the maximum number of edges of a simple graph
with $n$ vertices, $n$ even and large enough, and with exactly $k$ perfect
matchings is equal to $n^2/4+c_k$ \cite{DudekSchmitt2012}, where $c_k \leq k$ and
$c_k$ is positive for $k>1$ \cite{HartkeStoleeWestYancey2013}. Another ``extremal'' result of a similar flavour states that for every simple graph $G$ on $n$ vertices and
$m$ edges there exists a graph $H$ on $n$ vertices and $m$ edges with $|\mathfrak{M}(H)| \leq |\mathfrak{M}(G)|$ such that $H$ is a {\em threshold graph},
that is, it admits a clique $K$ such that the vertices from $V(H) \setminus K$ are independent and their neighborhoods form a chain with respect to $\subseteq$.
This has been used to determine the minimum number of perfect matchings in a simple graph on $n$ vertices and $m$ edges \cite{GrossKahlSaccoman2010};
although being a minimizing result at first glance, that number is trivially $0$ if $m \leq {n \choose 2} - (n-1)$, so that the interesting part of the analysis
is concerned with extremely dense graphs.
Among the few structural results on graphs with a fixed or even only a small number of perfect matchings let us mention
{\sc Lov\'asz}'s Cathedral Theorem (see Chapter 5 in \cite{LovaszPlummer1986}), which characterizes the {\em maximal} graphs having exactly $k$ perfect matchings, and
{\sc Kotzig}'s classic theorem that every connected graph with a unique perfect matching admits a bridge from that matching \cite{Kotzig1959}.
The latter theorem has been used recently to prove that a graph $G$ without three pairwise nonadjacent vertices and exactly one optimal coloring (in terms
of the chromatic number) has a shallow clique minor of order at least $|V(G)|/2$ \cite{Kriesell2016}, which supports {\sc Seymour}'s conjecture that
every graph $G$ without three pairwise nonadjacent vertices in general admits a shallow clique minor of order at least $|V(G)|/2$.
By getting more structural insight into graphs (and also hypergraphs) with only a few perfect matchings --- as provided by our main result ---
it may be possible to generalize the results from \cite{Kriesell2016}.

Let us close this section with two simple observations.
First note that every edge incident with some vertex $x$ in a minimally $k$-matchable graph
must be contained in at least one perfect matching; since every perfect matching contains exactly one edge incident with $x$, the degree of $x$, and, hence,
the maximum degree of $G$, is bounded from above by $|{\frak M}(G)|$.
As we have seen above, the number of perfect matchings of a minimally $k$-matchable graph $G$ can be larger than $k$, but the
following Lemma bounds it by $2k-2$ (and bounds, at the same time, the maximum degree $\Delta(G)$ of $G$ by $k$).

\begin{lemma}
  \label{L1}
  Let $x$ be a vertex of a minimally $k$-matchable graph $G$ with degree $d:=d_G(x) \geq 2$.
  Then $|{\frak M}(G)| \leq \frac{d}{d-1} \cdot (k-1)$. In particular, $|{\frak M}(G)| \leq 2k-2$ and $\Delta(G) \leq k$.
\end{lemma}

{\bf Proof.}
Let $x$ be a vertex of degree $d:=d_G(x) \geq 2$, and let $J$ be the set of edges of $G$ incident with $x$; so $|J|=d$.
For $e \in J$, let $m_e$ denote the number of perfect matchings from $|{\frak M}(G)|$ containing $e$.
Consequently, $|{\frak M}(G)|=\sum_{e \in J} m_e =:s$. Since $G-e$ has $s-m_e$ perfect matchings
and $G$ is minimally $k$-matchable, we get $s-m_e \leq k-1$. Taking the sum over all $e \in J$ on both sides
we get $d \cdot s - s \leq d \cdot (k-1)$, from which the statement of the Lemma follows.
Since $d/(d-1)$ is decreasing for increasing $d$, it is maximal for $d=2$, implying $|{\frak M}(G)| \leq 2k-2$.
Since $s \geq k$ by assumption to $G$ we derive $d_G(x) \leq k$ for all vertices and hence $\Delta(G) \leq k$.
\hspace*{\fill}$\Box$

The following Lemma implies easily the formally stronger version of Theorem \ref{T1} that
for every $k \geq 1$ there exists a finite set of graphs ${\frak G}_k$ such that a graph is minimally $k$-matchable
{\em if and only if} it is isomorphic to the disjoint union of an odd subdivision of some graph from ${\frak G}_k$ and any number of copies of $K_2$.

\begin{lemma}
  \label{L2}
  Let $G$ be the disjoint union of an odd subdivision of some graph $H$ and any number of copies of $K_2$.
  Then $|{\frak M}(G)|=|{\frak M}(H)|$.
\end{lemma}

{\bf Proof.}
Suppose that $G$ has been obtained from $H$ by disjointly adding a single copy of $K_2$, and let $e$ be the edge of that $K_2$.
One checks readily that $\varphi:{\frak M}(H) \rightarrow {\frak M}(G)$, $\varphi(M):=M \cup \{e\}$, is a bijection.
Suppose that $G$ has been obtained from $H$ by replacing an edge $wz$ by a path $wxyz$ of length $3$, where
$x,y$ are new vertices. For a perfect matching $M$ of $H$, define $\psi(M):=(M \setminus \{wz\}) \cup \{wx,yz\}$ if $wz \in M$ and
$\psi(M):=M \cup \{xy\}$ if $wz \not\in M$. In either case, $\psi(M)$ is a perfect matching of $G$, and $\psi:{\frak M}(H) \rightarrow {\frak M}(G)$
constitutes a bijection. Since any disjoint union of an odd subdivision of $H$ and any number of copies of $K_2$
can be obtained by subsequently disjointly adding single copies of $K_2$ or replacing edges by paths of length $3$ with new
internal vertices, the statement of the Lemma follows by induction.
\hspace*{\fill}$\Box$

\section{Proof of Theorem \ref{T1}}

For a path $P$ and vertices $a,b$ from $P$, let $aPb$ denote the subpath of $P$ connecting $a$ and $b$.
We apply this notion to some cycles as well; to this end, such a cycle $C$ comes with a fixed orientation,
and for vertices $a \not= b$ from $C$, $aCb$ is the subpath from $a$ to $b$ of $C$ following that orientation;
we also refer to $aCb$ as the {\em $a,b$-segment} along $C$. By $C^{-1}$, we denote the cycle $C$ with
the orientation opposite to the given one (so the $a,b$-segment along $C$ is the $b,a$-segment along $C^{-1}$).
$R:=P_1 \dots P_k$ denotes the union (concatenation) of the paths $P_1,\dots,P_k$. If the $P_j$ are described
as subpaths of larger paths or segments along cycles by their end vertices, say, $P_i=a_iQ_ib_i$,
and if $b_i=a_{i+1}$ then we list only one of $b_i,a_{i+1}$ in the description of $R$; for example, we write
$aPbQc$ instead of $aPbbQc$. In all cases, $R$ will be a path or a cycle.

Let $M$ be a perfect matching of a graph $G$. A cycle $C$ is {\em $M$-alternating} if $M \cap E(C)$ is a perfect matching of $C$.
If $C$ is $M$-alternating then the symmetric difference $(M \setminus E(C)) \cup (E(C) \setminus M)$ of $M$ and $E(C)$
is a perfect matching, too, and we call it the matching obtained from $M$ by {\em exchanging} along $C$.
If $N$ is another perfect matching then a path $P$ is called {\em $N,M$-alternating} if $N$ is a perfect matching of $P$
and $E(P) \setminus N \subseteq M$; that is, $P$ starts and ends with an edge of $N$ and if $f,g$ are consecutive on $P$
then $f \in N \wedge g \in M$ or $f \in M \wedge g \in N$.

{\bf Proof of Theorem \ref{T1}.} \\
We do induction on $k$. The statement is obviously true for $k=1$, take ${\frak G}_1=\emptyset$.
Let $G$ be a minimally $(k+1)$-matchable graph.
We may assume that $G$ is not an odd subdivision of some smaller graph,
and that no component of $G$ is isomorphic to $K_2$.
Since $\Delta(G) \leq k+1$ by Lemma \ref{L1}, it suffices to find an upper bound for $|V(G)|$ in terms of $k$.

$G$ contains a spanning minimally $k$-matchable subgraph $H$.
By induction, $H$ is the disjoint union of an odd subdivision of some graph from ${\frak G}_k$ and some number of copies of $K_2$.
Let $F:=E(G) \setminus E(H)$, and let ${\frak N}:={\frak M}(G) \setminus {\frak M}(H)$.
Since ${\frak G}_k$ is finite by induction, it suffices to bound the length of the subdivision paths in $H$ and (which is much easier) the number of copies of $K_2$ in terms of $k$ from above.
If $F$ is empty then this is obvious; $G$ is then one of the graphs from ${\frak G}_k$.
Hence it suffices to consider the case that $F \not= \emptyset$, implying ${\frak N} \not= \emptyset$.

{\bf Claim 1.} $F \subseteq \bigcap {\frak N}$. In particular, $F$ is a matching, and no perfect matching of $G$ contains at least one but not all edges of $F$.

Suppose, to the contrary, that there exist $e \in F$ and $N \in {\frak N}$ 
such that $e \not\in N$.
Since ${\frak M}(G-e) \supseteq {\frak M}(H) \cup \{N\}$, $G-e$ has $k+1$ perfect matchings, contradicting the minimality of $G$.
This proves Claim 1.

Now let $M \in {\frak M}(H)$, $N \in {\frak N}$, and consider an $M$-alternating cycle $C$ in $H$ with some fixed orientation;
we orient the edges of $M$ accordingly.

\begin{figure}
  \begin{center} 
    \scalebox{0.5}{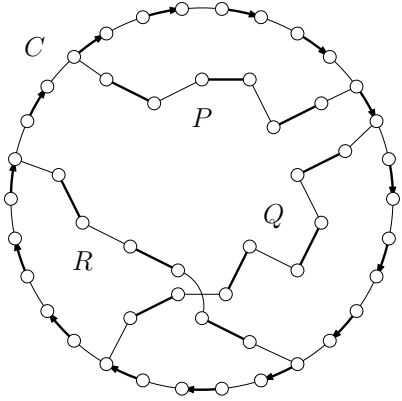} 
    \caption{\label{F1} An out-chord $P$, an in-chord $Q$, and an odd chord $R$ along the cycle $C$.
      Edges from the perfect matching $M$ are displayed fat and, on $C$, oriented according to the direction of $C$.
      Non-fat edges on the chords are necessarily from one and the same matching $N$, whereas non-fat edges on $C$
      may be anywhere outside $M$.
      $Q$ and $R$ cross, whereas $P,Q$ and $P,R$ do not cross.
      $R$ together with the lower segment along $C$ starting at $R$'s endvertices forms another $M$-alternating cycle.}
  \end{center}
\end{figure}

Deviant from standard notion, a {\em chord} of $C$ is an $N,M$-alternating (odd) path having only its end vertices in common with $C$.
Observe that for every edge $e \in N \setminus  E(C)$ incident with at least one vertex from $C$ there exists a chord starting with $e$.
Let $P$ be a chord, and let $a,b$ be its endvertices on $C$. Both $a,b$ are incident with a unique oriented edge $e,f$, respectively, from $M$.
If both $a,b$ are initial vertices of $e,f$, respectively, then we call $P$ an {\em out-chord}, if they are both terminal vertices then we
call $P$ an {\em in-chord}, and in the other cases $P$ is called an {\em odd chord}. $P$ is {\em external} if it contains at least one
edge from $F$ (which is then from $N$), and {\em internal} otherwise. $P$ {\em crosses} a chord $Q$ if the end vertices of
$Q$ are in distinct components of $C-\{a,b\}$; in that case, $Q$ crosses $P$, too. See Figure \ref{F1} for an example.

\DONOTTEX{
Observe that every chord is on a unique $N,M$-cycle in $G$.

{\bf Claim 2.} All external chords are on the same $N,M$-cycle.

For otherwise the perfect matching obtained from $M$ by exchangig $N,M$ on some $N,M$-cycle would contain at least
one but not all edges from $F$, so it would be a matching in ${\frak N}$ contradicting Claim 1. This proves Claim 2.
}
   
{\bf Claim 2.} If $P$ is an odd external chord then it is the only external chord.

$P$ can be extended to an $M$-alternating cycle by (exactly) one of the two paths connecting its end vertices in $C$.
The perfect matching obtained from $M$ by exchanging along this cycle 
would contain the edges from $E(P) \cap N$
but no other edges from $F$; it is from ${\frak N}$
(see, for example, the odd chord $R$ in Figure \ref{F1}), 
so that, by Claim 1, $F \subseteq E(P) \cap N$ follows;
in particular, there cannot be another external chord. This proves Claim 2.

{\bf Claim 3.} Suppose that some in-chord $P$ crosses some out-chord $Q$. 
Then either both $P,Q$ are internal, or there are no external chords distinct from $P,Q$.

\begin{figure}
  \begin{center} 
    \scalebox{0.5}{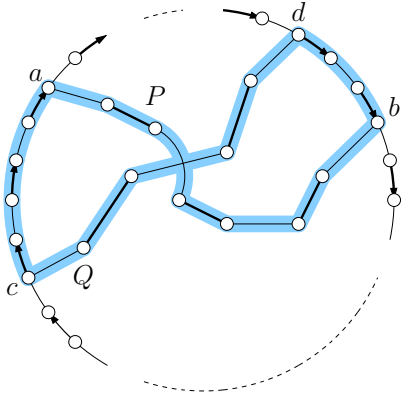} 
    \caption{\label{F2} An in-chord $P$ and an out-chord $Q$ which cross;
      their union with the linkage connecting the endvertices $a,c$ and the endvertices $b,d$
      forms another $M$-alternating cycle, underlayed in grey.}
  \end{center}
\end{figure}

$P \cup Q$ can be extended to an $M$-alternating cycle along $C$ by (exactly) one of the two linkages connecting their end vertices in $C$ (see Figure \ref{F2}).
The matching $M'$ obtained from $M$ by exchanging along this cycle would contain the edges from $(E(P) \cup E(Q)) \cap N$
but no other edges from $F$; if not both $P,Q$ are internal, then $M'$ is from ${\frak N}$, so that, by Claim 1, $F \subseteq (E(P) \cup E(Q)) \cap N$ follows;
in particular, there cannot be another external chord except for $P,Q$. This proves Claim 3.

We now turn to a more specific situation concerning $C$. Suppose that $D=x_0 x_1 \dots x_\ell$ is a subpath of $C$ 
of length $\ell \geq 6$ whose vertices have degree $2$ in $H$. We will show that if $\ell$ is large then we find a large number of $M$-alternating cycles
in $G$, each with an edge not in any of the others, from which we can construct a very large number of perfect matchings in $G$, contradicting Lemma \ref{L1}.

If $D$ contained an edge of $N$ then by Claim 1 both of its endvertices have degree $2$ in $G$, and from this it (easily) follows
that $G$ is an odd subdivision of a smaller graph, which has been excluded initially. Therefore, $D$ contains no edges from $N$.
Since $N$ is a perfect matching, every internal vertex $x_i$ of $D$ (that is: $i \in \{1,\dots,\ell-1\}$) is the end vertex of an external chord, say, $P_i$.
By Claim 2, these chords are in- or out-chords, and $P_i$ is an in-chord if and only if $P_{i+1}$ is an out-chord, for all $i \in \{1,\dots,\ell-2\}$.
By Claim 3, $P_i$ and $P_{i+1}$ do not cross, implying that $P_i$ and $P_j$ are distinct and do not cross, for all $i \not= j$ from $\{1,\dots,\ell-1\}$.
In particular, there are at least three external chords; Claim 2 thus implies that {\em there are no external odd chords} at all, and Claim 3 implies
that, in general, {\em an in-chord and an out-chord cannot cross unless they are both internal}.

{\bf Claim 4.}
No chord crosses three of the $P_i$.

\begin{figure}
  \begin{center} 
    \scalebox{0.5}{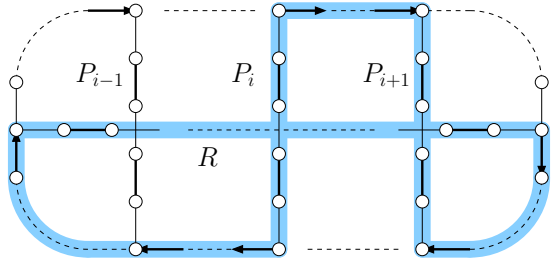} 
    \caption{\label{F3} Three consecutive chords $P_{i-1},P_i,P_{i+1}$ crossed by a (long, horizontal) odd chord $R$.
    $P_{i-1}$ and $P_{i+1}$ are in-chords, $P_i$ is an out-chord. The resulting $M$-alternating cycle is underlayed in grey.
    We get similar pictures if $P_{i-1},P_{i+1}$ were out-chords and $P_i$ was an in-chord, or if the edges from $M$ on $C$
    incident with the end vertices of $R$ were actually in the upper half.
    The picture does not determine the location of $x_{i-1},x_i,x_{i+1}$;
    if they are on the upper part of the picture then the $x_{i-1},x_i$-segment and the $x_i,x_{i+1}$-segment along $C$ each consist
    of a single edge only (whereas the picture suggests that these may be longer segments). 
    }
  \end{center}
\end{figure}

Suppose that some chord $R$ crosses three of the $P_i$. Then it crosses three consecutive of them, say $P_{i-1},P_i,P_{i+1}$.
Since at least one among them is an in-chord and at least one is an out-chord, $R$ must be an odd chord by Claim 3 and, thus, internal by Claim 2.
$R$ extends to an $M$-alternating cycle as follows:
We extend $R$ along its outgoing edge of $M$ along $C$ until we meet the first in-chord among $P_{i-1},P_i,P_{i+1}$,
follow that in-chord, exit it via its second in-edge on $C$,
follow $C$ opposite to its given orientation until we meet the next chord among $P_{i-1},P_i,P_{i+1}$, which is an out-chord,
traverse that out-chord, exit via its second out-edge on $C$, and close by traversing $C$ in its given orientation until we meet $R$
(see Figure \ref{F3} for an example).
By exchanging $M$ along this cycle we get a matching which contains the $N$-edges of two but not of all external chords,
violating Claim 1. This proves Claim 4.

Let $y_i$ denote the end vertex of $P_i$ distinct from $x_i$, and let $S_i$ denote the (closed)
$y_{i+1} y_i$-segment along $C$.

{\bf Claim 5.}
Let $i \in \{2,\ell-2\}$.
If $P_i$ is an out-chord then it is crossed by an internal odd chord or it is crossed by an out-chord with end vertices in $S_{i-1}$ and $S_i$.
If $P_i$ is an in-chord then it is crossed by an internal odd chord or it is crossed by an in-chord with end vertices in $S_{i-1}$ and $S_i$.

Suppose first that $P_i$ is an out-chord. The $y_i,x_i$-segment $D$ along $C$ has an odd number of vertices. An even number among them
is covered by edges from $N \cap E(C)$, so that an odd number among them is incident with an edge from $N$ not on $C$, i.~e. with an
end edge of some external chord. Since both $x_i,y_i$ are of the latter kind, there must be and odd number and, hence, at least one chord $Q$ starting in the
interior of $D$ and ending in $V(C) \setminus V(D)$, that is, $Q$ crosses $P_i$.
If $Q$ is odd then it is internal by Claim 2.
Otherwise, $Q$ must be an out-chord by Claim 3.
Again by Claim 3, $Q$ cannot cross the external in-chords $P_{i-1}$ or $P_{i+1}$,
so that its end vertices are in $S_{i-1}$ and $S_i$. This proves the first part of Claim 5, and, symmetrically, the second part follows.

{\bf Claim 6.}
Suppose that $x_i x_{i+1}$ is in $M$. Then there exists an $M$-alternating cycle distinct from $C$
in the subgraph $H_i$ formed by $C$ and all chords with both end vertices in $S_i \cup \dots \cup S_{i+5}$.

Observe that $x_{i+2} x_{i+3} \in M$ by construction. If $P_{i+2}$ or $P_{i+3}$ is crossed by an internal odd chord $S$ then its end vertices
are in $S_i \cup S_{i+1} \cup S_{i+2} \cup S_{i+3} \cup S_{i+4}$ by Claim 5; hence the unique $M$-alternating cycle in $C \cup S$ 
containing $S$ verifies Claim 6 in this case.
Hence we may suppose that neither $P_{i+2}$ nor $P_{i+3}$ is crossed by an internal odd chord.
If there was an internal odd chord $S$ with some end vertex in $S_{i+2}$ then its other end vertex would be in $S_{i+2}$, too,
and the unique $M$-alternating cycle in $C \cup S$ containing $S$ verified Claim 6 again.
Hence 
\begin{quote}
  all chords with end vertices in $S_{i+2}$ are in-chords or out-chords. $(\ast)$
\end{quote}

Suppose that there is an in- or out-chord $Q$ with both end vertices in $S_{i+2}$. Take it in such a way that the distance of its end vertices 
is as small as possible in the graph $S_{i+2}$.
Let  $a$ and $b$ be the end vertices of $Q$. Exactly one of $a,b$ is incident with an edge from
$M \cap E(aS_{i+2}b)$. Without loss of generality, let it be $a$; there is an $M,N$-alternating subpath of $aS_{i+2}b$ starting with $a$,
and we take a maximal one, say $S$; its end vertex $c$ distinct from $a$ is an internal vertex of $aS_{i+2}b$, and by maximality of $S$
the edge $e$ from $N$ incident with $c$ is not in $E(C)$; observe that $e \not= bc$ since the edge from $N$ incident with 
% gf $c$
$b$ %gf 
is on $Q$.
Hence there is a chord $R$ with end vertex $c$, and $R \not= Q$. 
It must either be an in-chord or an out-chord by 
$(\ast)$ %gf
as $c \in S_{i+2}$,
and since $S$ is an $M,N$-alternating path, we know that $R$ is an in-chord if $Q$ is an out-chord and $R$ is an out-chord if $Q$ is an in-chord.
By choice of $Q$, the end vertex $d$ of $R$ distinct from $c$ is not in $aS_{i+2}b$, so that $Q,R$ cross.
If $Q$ is an in-chord then 
% gf $bQaScRdCb$ 
$bQaScRdC^{-1}b$ %gf
is the desired $M$-alternating cycle:
In that case, $R$ is an out-chord, so it cannot cross the external in-chords $P_{i+3}$ and $P_{i+1}$, implying $d \in S_{i+1} \cup S_{i+2}$),
and both of $Q,R$ are internal by Claim 3.
If, otherwise, $Q$ is an out-chord then, symmetrically, $bQaScRdC^{-1}b$ is the desired $M$-alternating cycle.

\begin{figure}
  \begin{center} 
    \scalebox{0.4}{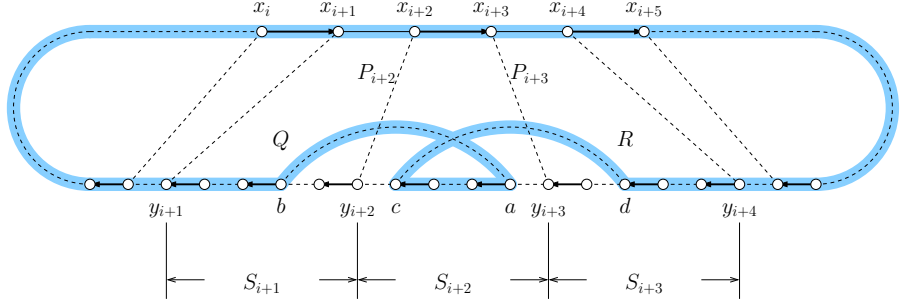} 
    \caption{\label{F4} Finding the desired $M$-alternating cycle in Claim 6 (underlayed in grey).
      Edges from $M \cap E(C)$ are displayed fat as before, dashed connections resemble paths of odd length.
      The vertices $x_i,\dots,x_{i+5}$ are consecutive on $C$, so there is ``no space'' for the end vertices
      of chords other than $P_i,\dots,P_{i+5}$ ``in between'' them. Some labels are omitted.}
  \end{center}
\end{figure}

Hence all chords with some end vertex in $S_{i+2}$ must cross $P_{i+2}$ or $P_{i+3}$.
$P_{i+2}$ is an out-chord, so that, by Claim 5 and $(\ast)$, it is crossed by an out-chord $Q$ with end vertices $b \in S_{i+1}$ and $a \in S_{i+2}$;
$a$ is adjacent with an edge from $M \cap E(aS_{i+2}b)$. As in the previous paragraph, there exists a maximal $M,N$-alternating path 
in $aS_{i+2}b$ starting with $a$ and ending with a vertex $c \not=a$. Since $y_{i+2}$ is end vertex of an out-chord, we see
that $c$ is an inner vertex of $aCy_{i+2}$. As above, there is an in-chord $R$ with end vertex $c$. $R$ crosses either $P_{i+2}$ or $P_{i+3}$,
but it cannot cross the external out-chord $P_{i+2}$, so that it must cross $P_{i+3}$.
But then the end vertex $d$ of $R$ distinct from $c$ is
in $S_{i+3}$ as $R$ cannot cross the external out-chord $P_{i+2}$. It follows that $Q,R$ cross, so they are internal by Claim 3,
and $bQaScRdC^{-1}b$ (where $d$ is the end vertex of $R$ distinct from $c$) is the desired cycle. Figure \ref{F4} illustrates the process.
This proves Claim 6.

\DONOTTEX{
If $b$ is closer to $y_{i+2}$ than $c$ then $P,Q$ do not cross and
$x_i x_{i+1} P_{i+1} y_{i+1} C^{-1} d R c C b Q a C^{-1} y_i P_i x_i$ is an $M$-alternating cycle;
the matching obtained from $M$ by exchanging along that cycle produces a matching which contains all edges from $F$ in $P_i,P_{i+1}$
but no edges from $F$ in any other $P_j$, contradicting Claim 1.

Hence $c$ is closer to $y_i$ than $b$, implying that $P,Q$ cross;
in that case, $aQ b C^{-1} d R b C a$ is an $M$-alternating cycle, proving Claim 7.}

Now consider an {\em arbitrary} path $x_0,\dots,x_\ell$ of vertices of degree $2$ in $G$ and observe that it is contained in some
$M$-alternating cycle $C$, to which we apply the considerations following Claim 1.
We construct an upper bound for $\ell$ in terms of $k$.
There exists a $d$ such that $\ell-1 \geq 6d+1$ but $\ell-1 <6(d+1)+1$.
Then for some $j_0 \in \{1,2\}$, $x_{j_0} x_{j_0+1}$ is in $M$. For $j \in \{0,\dots,d-1\}$ and $i:=j_0+6 \cdot j$
there exists an $M$-alternating cycle $C_j$ in $H_i$ as in Claim 7, and the sets $E(C_j) \setminus E(C)$ are nonempty and pairwise disjoint. 
For every $J \subseteq \{0,\dots,d-1\}$, let $M_J$ be the symmetric difference of $M$ and $(C_j)_{j \in J}$, that is
$M_J:=\{e \in E(G):$ $e$ is contained in an odd number of $M,(C_j)_{j \in J}\}$, is a perfect matching of $H$, and $M_J \not= M_{J'}$ for
$J \not= J'$. By Lemma \ref{L1} and Lemma \ref{L2}, $H$ has at most $2k-2$ perfect matchings,
so that $2^d \leq 2k-2$, that is, $d \leq \log_2 (k-1)+1$.
It follows $\ell \leq 6(d+1) \leq 6 \log_2(k-1)+12$. 

Recall that $H$ is the disjoint union of an odd subdivision of some graph from ${\frak G}_k$, say, $H_0$, and some number, say $q$, of copies of $K_2$.
Suppose that $e$ is the edge of one of the latter copies of $K_2$, and let us assume, to the contrary, that $e$ had no parallel edges in $G$.
If one of the endvertices had degree $1$ in $G$ then $e$ would be contained in every perfect matching of $G$;
if there was an edge $f \not=e$ incident with $e$ then it cannot be contained in any perfect matching of $G$,
so that ${\frak M}(G-f)={\frak M}(G)$, contradiction; therefore, both endvertices in $G$ had degree $1$, contradicting the initial assumption that $G$
has no components isomorphic to $K_2$. Consequently, both end vertices of $e$ were incident with edges from $F$, which is a matching by Claim 1.
By assumption to $e$, these edges were distinct, and both end vertices of $e$ had degree $2$ in $G$; from this it easily follows
that $G$ is an odd subdivision of some smaller graph, contradiction.
Therefore, every edge forming a copy of $K_2$ in $H$ must have at least one --- and, hence, exactly one --- parallel in $G$,
so that $G$ had at least $2^q$ perfect matchings; it follows that $2^q \leq 2k$ by Lemma \ref{L1}, implying $q \leq \log_2 k +1$.

Since every edge in $H_0$ is subdivided by at most $6 \log_2(k-1)+12$ vertices,
$|V(G)| \leq |V(H_0)|+(6 \log_2(k-1)+12) \cdot |E(H_0)| + 2 \log_2 k +2$. 
As ${\frak G}_k$ is finite, we get
$|V(G)| \leq f(k)$ with
$f(k):=\max\{|V(H')|+(6 \log_2(k-1)+12) \cdot |E(H')| + 2 \log_2 k +2:\,H' \in {\frak G}_k\}$.
(And, as already mentioned above, $|E(G)| \leq (k+1) \cdot f(k)$ by Lemma \ref{L1}.)
\hspace*{\fill}$\Box$

\section{Minimally \boldmath$3$-matchable graphs}

Let us finish by describing the set ${\frak G}_3$ by specializing (and, thus, partly illustrating) the ideas of the proof of Theorem \ref{T1}.
We may assume that $G \in {\mathfrak G}_3$ is minimally $3$-matchable, i.e.\ 
\begin{quote}
$G$ is not an odd subdivision of some smaller graph, $(\dagger)$
\end{quote}
and that no component of $G$ is isomorphic to $K_2$.
$G$ contains a spanning minimally $2$-matchable subgraph $H$, and, according to Claim 1, $F:=E(G)-E(H)$ is a (not necessarily perfect) matching.
As ${\frak G}_2=\{C_2\}$, $H$ is the disjoint union of an even cycle $H_0$ and $q$ copies of $K_2$.
By repeating the arguments in the end of the proof of Theorem \ref{T1} we see that
any of these copies must have exactly one parallel edge in $G$. 
As $H_0$ has two matchings we see that $q \leq 1$, for otherwise $G$ had at least $8$ matchings, contradicting Lemma \ref{L1}.
Moreover, if $q=1$ then there is no edge $e \in F$ connecting two vertices from $H_0$, for otherwise $G-e$ contained two disjoint even cycles and thus still had four matchings;
in that case we deduce that $G$ consists of two disjoint $2$-cycles.

If, otherwise, $q=0$, then $H=H_0$, and there must be at least one edge $e \in F$ connecting two distinct vertices from $H$ in $G$ (as $H$ has only two perfect matchings).
In order to apply the chord notion of the previous section, let us take a matching $M$ of $H$ and another matching $N$ not from $H$, and fix an orientation of 
$H \;(= C)$. %gf  
If $e$ forms and external odd chord then $H+e$ already contains three disjoint matchings, so that $H+e=G$, and, by 
$(\dagger)$, %gf 
$G$ is the graph on two vertices with three parallel edges,
sometimes called the {\em theta graph} (see also Claim 2 above). So we may assume without loss of generality that all edges from $F$ constitute external in- or out-chords.
We take $e=xy$ such that the length of the $x,y$-segment along $C=H$ is minimized. Without loss of generality, $e$ is an out-chord (otherwise we reverse the orientations and $x,y$).
There exists a maximal subpath $P$ of $S$ of starting at $x$ with the out-edge from $M$ and then alternately using edges from $N$ and $M$.
Its endvertex $z$ distinct from $x$ is an interior vertex of $S$ and its final edge is again from $M$.
Hence there exists an edge $g$ from $F$ constituting an in-chord. By choice of $e$, $g$ crosses $e$.
Now it follows (as in Claim 3), that there are no further external chords at all. Consequently, by $(\dagger)$, $G$ is the complete graph $K_4$.
Hence we proved:

\begin{theorem}
  \label{T2}
  Every minimally $3$-matchable graph is isomorphic to the disjoint union of any number of copies of $K_2$ and
  either two $2$-cycles, or an odd subdivision of the theta graph, or an odd subdivision $K_4$.
\end{theorem}

It is possible to restate Theorem \ref{T1} and its specializations in terms of chambers as used in connection with {\sc Lov\'asz}'s Cathedral Theorem (see Chapter 5 in \cite{LovaszPlummer1986}).
Let us do this for Theorem \ref{T2}. According to \cite{HartkeStoleeWestYancey2013}, a {\em chamber} is the vertex set of a connected component of the spanning subgraph
$H:=(V(G),\bigcup \mathfrak{M}(G))$ formed by all edges of perfect matchings.
Now if $G$ is a graph with exactly three perfect matchings we know that $H$ is minimally $3$-matchable, so that, apart from chambers spanned by edges
in all three matchings, $G$ has either two further chambers spanned by an even cycle each, or a single further chamber spanned by a totally odd subdivision of the theta graph,
or a single further chamber spanned by an odd subdivision of $K_4$. Analogously, one could think of Theorem \ref{T1} as a classification theorem for graphs with exactly $k$ perfect matchings.

{\bf Addresses of the authors.}

\nopagebreak

\parbox{7cm}{
{\sc Ga\v{s}per Fijav\v{z}} \\
Faculty of Computer and Information Science\\
University of Ljubljana \\
Ve\v{c}na pot 113 \\
1000 Ljubljana \\
Slovenia}
\hspace*{\fill}
\parbox{48mm}{
{\sc Matthias Kriesell} \\
Institut f\"ur Mathematik\\
Technische Universit\"at Ilmenau \\
Weimarer Stra{\ss}e 25 \\
98693 Ilmenau \\
Germany}

\end{document}